\documentclass[twoside,centertags]{amsart}
\usepackage[margin=0.8in]{geometry}
\newtheorem{innercustomthm}{Theorem}
\newenvironment{customthm}[1]
  {\renewcommand\theinnercustomthm{#1}\innercustomthm}
  {\endinnercustomthm}
  
  \usepackage{fancyhdr}            
\fancypagestyle{plain}{
    \lhead{}
    \fancyhead[R]{\thepage}
    \fancyhead[L]{}
    
    \fancyfoot{}
}

\usepackage[cmtip,color,all]{xy}
\usepackage{amsmath}
\usepackage{amsfonts}
\usepackage{amssymb}
\usepackage{amsthm}
\usepackage{graphicx}
\usepackage{mathrsfs}
\usepackage{hyperref}
\usepackage{tikz-cd}
\usepackage{xcolor}
\usepackage{bbold}
\usepackage{pifont}

\theoremstyle{definition}

\theoremstyle{definition}

\numberwithin{equation}{section}

\begin{document}
\begin{abstract}
 We explain how the undecidability of the quasi-isomorphism and derived Morita equivalence problems for semi-free finitely presented dg associative algebras follows as an immediate consequence of an extension of Adams' cobar theorem to non-simply connected spaces and the triviality problem for finitely presented groups.
\end{abstract}
\title[Undecidability of the quasi-isomorphism problem for semi-free dg algebras via the cobar construction]{Undecidability of the quasi-isomorphism problem for semi-free dg algebras via the cobar construction}

\author[M. Rivera]{Manuel Rivera}
\address{\newline M.R., Department of Mathematics, Purdue University}
\email{\href{mailto:manuelr@purdue.edu}{manuelr@purdue.edu}}

\maketitle
In \cite{MR}, a proof to Theorem \ref{thm2} below is presented based on autonomous responses of \textit{Aletheia}, an AI math research agent, providing a solution to the first half of Problem 5.16 from the K3 Problem List in Low-Dimensional Topology \cite{K}. In this note, I explain how Theorem \ref{thm2} is a straightforward consequence of the main result of \cite{RiZe16} (Theorem \ref{thm1} below) and the undecidability of the triviality problem for finitely presented groups. I believe the argument is conceptually clarifying. In fact, given any finitely presented group, the cobar construction algorithmically provides a semi-free finitely presented differential graded (dg) algebra whose zeroth Hochschild homology is the free module on the set of conjugacy classes of the group. The \textit{Aletheia} raw output linked in \cite{MR} contains an attempt to a solution using \cite{RiZe16}, which the authors do not discuss. We first fix notation and recall the main theorem of \cite{RiZe16}. Fix a (non-trivial) commutative ring with unit $R$. 
\begin{enumerate}
    \item Denote by $\mathbb{\Omega} \colon \mathsf{dgCoalg}_R \to \mathsf{dgAlg}_R$
the \textit{cobar} functor from the category of dg coassociative coaugmented $R$-coalgebras to the category of dg associative augmented $R$-algebras. For any $C \in \mathsf{dgCoalg}_R$, the dg algebra $\mathbb{\Omega}(C)$ is \textit{semi-free}, its underlying graded algebra being the free associative graded algebra on $s^{-1}\overline{C}$, where $\overline{C}$ is the cokernel of the coaugmentation and $s^{-1}$ the ``shift down by $1$" functor. The differential is defined by extending the differential and coproduct of $C$ as a graded derivation.

\item For any simplicial set $X$, denote by $C_*(X;R)$ the dg $R$-coalgebra of normalized chains on $X$ with Alexander-Whitney coproduct. A choice of vertex in $X$ gives a coaugmentation $R \to C_*(X;R)$. For any topological space $Y$, denote by $S_*(Y;R)$ the dg $R$-module of normalized singular chains on $Y$. A topological monoid structure on $Y$ induces a dg algebra structure on $S_*(Y,R)$. 

\item If $X$ is a simplicial set with a single vertex and $\sigma \in X_{n>0}$ a non-degenerate $n$-simplex, the differential $d \colon \mathbb{\Omega}(C_*(X;R)) \to \mathbb{\Omega}(C_*(X;R))$
is given explicitly by 
\[ d (s^{-1}\sigma) = \sum_{i=0}^n(-1)^{i+1}s^{-1}\sigma(0,...,\hat{i},...,n) + \sum_{j=1}^{n-1} (-1)^j s^{-1}\sigma(0,...,j) \otimes s^{-1} \sigma(j,...,n),\]
where $\sigma(i_0,...,i_k)$ denotes the $k$-simplex obtained by restricting $\sigma$ to the vertices $i_0,...,i_k$.
\item Denote by $ \Omega \colon \mathsf{Top}^* \to \mathsf{Mon}_{\mathsf{Top}}$ the \textit{based (Moore) loops} functor from the category of pointed topological spaces to the category of topological monoids. 

\item The \textit{homotopy category} of a simplicial set $X$ is the category $\pi(X)$ obtained by applying the left adjoint of the nerve functor to $X$.

\item Two dg associative algebras are \textit{quasi-isomorphic} if they are connected by a zig-zag of maps of dg associative algebras each inducing an isomorphism on homology. 
\end{enumerate}
\begin{customthm}{1}[\cite{RiZe16}]\label{thm1} 
If $X$ is a simplicial set with a single vertex and $\pi(X)$ is a group, then the dg associative algebras
$\mathbb{\Omega}C_*(X;R)$ and $S_*(\Omega|X|;R)$ are naturally quasi-isomorphic.
\end{customthm}
\noindent The above theorem was originally stated for Kan complexes in \cite{RiZe16}. However, the composition of functors $\mathbb{\Omega} \circ C_*$ sends \textit{categorical equivalences} (also known as \textit{Joyal equivalences}) between simplicial sets to quasi-isomorphisms of dg algebras \cite[Proposition 8.2]{RiZe16}.  Any simplicial set whose homotopy category is a groupoid is categorically equivalent to a Kan complex. The notion of categorical equivalence is strictly stronger than weak homotopy equivalence. We now deduce the main theorem. 
\begin{customthm}{2}\label{thm2} If $R$ is Turing computable, the problem of determining whether two semi-free finitely presentable differential graded associative $R$-algebras are quasi-isomorphic is undecidable.
\end{customthm}
\noindent \textit{Proof of Theorem 2.}
We reduce the statement to the triviality problem for finitely presented groups, which is known to be undecidable \cite{Ra58}. Given any finitely presented group $(G,P)$ we associate two explicit semi-free, finitely presented, dg algebras $A_P$ and $B_P$ which are quasi-isomorphic if and only if $G$ is trivial. First, one can effectively associate to $(G,P)$ a $2$-skeletal simplicial set $K'_P$ with a single vertex and a finite number of simplices modeling the presentation complex of $P$. Glue a copy of $\Delta^3 / \sim$, where $\sim$ collapses $[0,2]$ to a vertex $[0]=[2]$ and $[1,3]$ to $[1]=[3]$, along the $1$-simplex $[1,2]$ to each $1$-simplex of $K'_P$. Denote by $K_P$ the resulting simplicial set. Then $\pi(K_P)\cong \pi_1(|K_P|)\cong \pi_1(|K'_P|)\cong G$. Define $A_P=\mathbb{\Omega} C_*(K_P;R)$. By Theorem \ref{thm1}, $A_P$ is quasi-isomorphic to $A'_P=S_*(\Omega|K_P|;R)$. If $G$ is the trivial group, then $A_P$ is quasi-isomorphic to $B_P=\mathbb{\Omega} C_*(\bigvee_{i=1}^{k(P)}S^2;R)$, where $k(P)$ the number of relations minus the number of generators in $P$ and $S^2$ is the simplicial set $\Delta^2/\partial\Delta^2$. Conversely,  if $A_P$  and $B_P$ are quasi-isomorphic, by Theorem \ref{thm1} and the quasi-isomorphism invariance of Hochschild homology, we have isomorphisms
\[HH_*(A'_P, A'_P) \cong HH_*(A_P,A_P)  \cong HH_*(B_P, B_P) \cong HH_*(B'_P,B'_P),\]
where $B'_P=  S_*(\Omega|\bigvee_{i=1}^{k(P)}S^2|;R)$. By section V.1 of \cite{G85}, there are isomorphisms
\[H_*(L|K_P|;R) \cong HH_*(A'_P, A'_P) \cong HH_*(B'_P,B'_P) \cong H_*(L|\bigvee_{i=1}^{k(P)}S^2|;R) \]
where $L$ denotes the free loop space functor. Since the zeroth homology of the free loop space is the free $R$-module generated by the set of conjugacy classes of the fundamental group, it follows that $G=\pi_1(|K_P|)$ has a single conjugacy class, so $G$ is trivial. \qed
\\
\\
\textbf{Remarks.}

\begin{enumerate}
    
\item If two dg algebras are quasi-isomorphic then they are derived Morita equivalent and Hochschild homology is derived Morita equivalence invariant. Hence, Theorem \ref{thm2} holds when ``quasi-isomorphic" is replaced by ``derived Morita equivalent".

\item The cobar functor $\mathbb{\Omega}$ does not send quasi-isomorphisms of dg coalgebras to quasi-isomorphisms of dg algebras. If $X$ is a simplicial set with exactly one vertex, one non-degenerate $1$-simplex, and every other simplex degenerate, the natural map $f \colon C_*(X;R) \to S_*(|X|;R)$ is a quasi-isomorphism of dg coaugmented coalgebras but $H_0(\mathbb{\Omega}(f))$ is the inclusion $R[x] \to R[x,x^{-1}]$.
\item If $C$ and $C'$ are \textit{simply connected} dg coalgebras (i.e. $R=C_0=C'_0$, $0=C_1=C'_1$, and $0=C_i=C'_i$ for $i<0$), any quasi-isomorphism $f \colon C \to C'$ induces  a quasi-isomorphism 
$\mathbb{\Omega}(f)\colon \mathbb{\Omega}(C) \to \mathbb{\Omega}(C')$. This follows from a standard spectral sequence argument.

\item A quasi-isomorphism of \textit{homologically} simply connected dg coaugmented coalgebras may not induce a quasi-isomorphism of dg algebras after applying $\mathbb{\Omega}$. Let $Y$ be a pointed space with non-trivial perfect fundamental group and let $Y^+$ be its Quillen plus construction. Then $0=H_1(Y;R)=H_1(Y^{+};R)$ and there is a quasi-isomorphism $S_*(Y;R) \to S_*(Y^{+};R)$, that does not induce a quasi-isomorphism of dg algebras after applying $\mathbb{\Omega}$.
\item Two better questions, which remain open at moment, are the following: \textit{Suppose that $\mathcal{A}$ and $\mathcal{B}$ are semi-free finitely presented dg associative algebras that are connected, i.e. $\mathcal{A}_0=R=\mathcal{B}_0$, and non-negatively graded. Is the problem of determining whether $\mathcal{A}$ and $\mathcal{B}$ are quasi-isomorphic/derived Morita equivalent algorithmically decidable? What if $\mathcal{A}$ and $\mathcal{B}$ are required to be semi-free finitely presented dg \emph{commutative} algebras?} See \cite{A85} for related results.
 \end{enumerate}

\noindent \textbf{Acknowledgments.} The author would like to thank Eric Samperton for fruitful discussions about decision problems.


\begin{thebibliography}{10}


\bibitem{A85} D. ~Anick,
\emph{Diophantine equations, Hilbert series, and undecidable spaces.} Ann. of Math. (2) 122 (1985), no. 1, 87--112.


\bibitem{G85} T. ~Goodwillie, \emph{Cyclic homology, derivations, and the free loopspace.} Topology 24 (1985), no. 2, 187--215.

\bibitem{K} R.I.~Baykur, R. C. ~Kirby, and D. ~Ruberman, editors, \emph{K3: A New Problem List in Low-Dimensional Topology},
volume 295 of Mathematical Surveys and Monographs, American Mathematical Society, Providence, RI, (2026).

\bibitem{MR} C. ~Manolescu, N. Rozenblyum, \emph{Undecidability problems for semi-free DG algebras}, arXiv:2605.08122, (2026).



\bibitem{RiZe16}
M.~Rivera and M.~Zeinalian, \emph{Cubical rigidification, the cobar construction and the based loop space.} Algebr. Geom. Topol. 18 (2018), no. 7, 3789--3820.

\bibitem{Ra58}
M.O.~Rabin, \emph{Recursive unsolvability of group theoretic problems.} Ann. of Math. (2) 67 (1958), 172--194.

\end{thebibliography}
\end{document}